%% file: TNOGOFknots-final.tex
\documentclass{amsart}
\usepackage{amsmath,amsthm}
\usepackage{graphics}
\usepackage{color}
\usepackage{epsfig}

\graphicspath{{./Figures/}}

\theoremstyle{plain}
\newtheorem{thm}{Theorem}[section]

\newtheorem{prop}[thm]{Proposition}
\newtheorem{lemma}[thm]{Lemma}
\newtheorem{claim}[thm]{Claim}

\newtheorem*{thmmain}{Theorem~\ref{thm:main}}
\newtheorem*{main2}{Theorem~\ref{thm:main2}}
\newtheorem*{main3}{Theorem~\ref{thm:main3}}

\theoremstyle{definition}
\newtheorem{remark}[thm]{Remark}

\newcommand{\B}{\ensuremath{{\mathfrak b}}}
\newcommand{\comment}[1]{}

\newcommand{\bdry}{\ensuremath{\partial}}

\newcommand{\Z}{\ensuremath{\mathbb{Z}}}

\begin{document}
\author{Kenneth L. Baker}
\author{Jesse E. Johnson}
\author{Elizabeth A. Klodginski}
\address{Department of Mathematics\\ University of Georgia\\ Athens GA 30602}
\email{kb@math.uga.edu}
\address{Department of Mathematics\\ University of California\\ Davis, CA 95616}
\email{jjohnson@math.ucdavis.edu}
\address{Department of Mathematics\\ University of California\\ Davis, CA 95616}
\email{eklodginski@math.ucdavis.edu}
\title{Tunnel number one, genus one fibered knots
}

\thanks{The first author was supported by NSF VIGRE grant DMS-0089927.  The second and third authors were supported by NSF VIGRE grant DMS-0135345.}



\begin{abstract}
We determine the genus one fibered knots in lens spaces that
have tunnel number one.  We also show that every tunnel number one, once-punctured torus bundle is the result of Dehn filling a component of the Whitehead link in the 3-sphere.
\end{abstract}

\maketitle

\section{Introduction}
A null homologous knot $K$ in a $3$--manifold $M$ is a {\em genus
one fibered knot}, GOF-knot for short, if $M-N(K)$ is a
once-punctured torus bundle over the circle whose monodromy is the
identity on the boundary of the fiber and $K$ is ambient isotopic in
$M$ to the boundary of a fiber.

We say the knot $K$ in $M$ has {\em tunnel number one} if there is a
properly embedded arc $\tau$ in $M-N(K)$ such that $M-N(K)-N(\tau)$
is a genus $2$ handlebody.  An arc such as $\tau$ is called an {\em
unknotting tunnel}.  Similarly, a manifold with toroidal boundary is
{\em tunnel number one} if it admits a genus-two Heegaard splitting.
Thus a knot is tunnel number one if and only if its complement is
tunnel number one.

In the genus $1$ Heegaard surface of $L(p,1)$, $p \neq 1$, there is
a unique link that bounds an annulus in each solid torus.  This two-component link
is called the {\em $p$--Hopf link} and is fibered with
monodromy $p$ Dehn twists along the core curve of the fiber.    
We refer to the fiber of a $p$--Hopf link as a {\em $p$--Hopf band}.  
In this terminology, the standard positive and negative Hopf bands
in $S^3 = L(+1, 1)=L(-1,1)$ are the $(+1)$-- and $(-1)$--Hopf bands,
respectively.

Gonz\'alez-Acu\~na~\cite{ga:dcok} shows that the trefoil (and its
mirror) and the figure eight knot are the only GOF-knots in $S^3$.
These knots arise as the boundary of the plumbing of two
$(\pm1)$--Hopf bands.  For both knots, a transverse arc on one of the
plumbed Hopf bands is an unknotting tunnel, so both have tunnel
number one.

GOF-knots in lens spaces were first studied by Morimoto in \cite{morimoto:gofkils}, and were classified by the first author of the present article in \cite{baker:cgofkils}.  Unlike in $S^3$, these knots do not all have tunnel number one.  In particular we prove the following:

\begin{thm}\label{thm:main} 
Every tunnel number one, genus one fibered knot in a lens space is
the plumbing of an $r$--Hopf band and a $(\pm 1)$--Hopf band which is contained in the
lens space $L(r,1)$, with one exception.  Up to mirroring, this
exception is the genus one fibered knot in $L(7,2)$ that arises as
$(-1)$--Dehn surgery on the plumbing of a $7$--Hopf band and a
$(+1)$--Hopf band.
\end{thm}

Every GOF-knot has tunnel number 1 or 2.  Thus Theorem~\ref{thm:main} determines the tunnel number of every GOF-knot in a lens space.

We also obtain
\begin{thm}\label{thm:main2}
Every tunnel number one, once-punctured torus bundle is the complement of
a GOF-knot in a lens space of the form $L(r,1)$.
\end{thm}

\begin{thm}\label{thm:main3}
Every tunnel number one, once-punctured torus bundle is the $(r/1)$--Dehn filling of a boundary component of the exterior of the Whitehead link for some integer $r$.
\end{thm}

This article begins by determining the monodromy of a tunnel number
one, once-punctured torus bundle in Section~\ref{sec:2}.  We then
translate these monodromies into the language of closed $3$--braids
in Section~\ref{sec:3}.  In particular we determine which
$3$--braids in a solid torus have double branched covers producing
tunnel number one, once-punctured torus bundles.
  In Section~\ref{sec:4} we discuss the presentations of two bridge links as closed $3$--braids in $S^3$.  In Section~\ref{sec:5} we determine which of the braids from Section~\ref{sec:3} are two bridge links.  The proof of Theorem~\ref{thm:main} is presented in Section~\ref{sec:6}, and Theorems~\ref{thm:main2} and \ref{thm:main3} are proved in Section~\ref{sec:7}.

We refer the reader to \cite{bz:knots} for background regarding
fibered knots, braids,  two bridge links, lens spaces, and double
coverings of $S^3$ branched over a link.  
Also, recall that plumbing is a local operation generalizing the connect sum.  
In particular, the plumbing of a fibered link in $S^3$ with a fibered link in another $3$--manifold $M$ produces a new fibered link in $M$.  
See for example Gabai's geometric description of the yet more generalized Murasugi sum in \cite{gabai:msiango}. 

The authors would like to thank Alan Reid for a helpful conversation
regarding the Whitehead link.

\section{The monodromy of a tunnel number one, once-punctured torus bundle.}\label{sec:2}

Heegaard splittings of closed torus bundles over the circle were studied in
detail by Cooper and Scharlemann~\cite{cooper-scharlemann}.  In
particular, their work characterized genus two Heegaard splittings of such
bundles.  Their method transfers almost directly to the case of
once-punctured torus bundles which are tunnel number one.

\begin{lemma}
\label{monodromylem} 
Let $M$ be a once-punctured torus bundle over the circle with once-punctured torus fiber $T$ and monodomy $\phi$.  Further assume $M$ allows a genus two Heegaard splitting.  
Then there is a pair of simple closed curves $\alpha_1$,
$\alpha_2$ in $T$ such that $\alpha_1 \cap \alpha_2$ is a single
point and $\phi$ sends $\alpha_1$ onto $\alpha_2$.  The map $\phi$
is isotopic to $(s_2 s_1 s_2)^{\pm 1} s_1^n$ where $s_1$ is a Dehn
twist along $\alpha_1$ and $s_2$ is a Dehn twist along $\alpha_2$.
\end{lemma}

\begin{proof}
We will show that the proof for closed torus bundles in Theorem~4.2
of \cite{cooper-scharlemann} works equally well in the once-punctured
case.  Because this method has been described in detail elsewhere,
we will give only an outline of the setup and leave many of the
details to the reader.  A similar exposition for general surface bundles
can also be found in \cite{bachman-schleimer}.

Let $(\Sigma, H_1, H_2)$ be a genus-two Heegaard splitting for $M$.
Assume $H_1$ is a compression body and $H_2$ a handlebody. A spine
$K_1$ for $H_1$ consists of $\partial_- H_1$ and an arc properly embedded in $H_1$ such that
the complement in $H_1$ of $K_1$ is homeomorphic to $\partial_+ H_1
\times (-1,0]$.  A spine $K_2$ for $H_2$ is a graph whose complement
is homeomorphic to $\partial H_2 \times [0,1)$. The Heegaard
splitting determines a continuous one-parameter family of embedded, pairwise disjoint
surfaces $\{\Sigma_x : x \in (-1,1)\}$ such that as $x$ approaches $-1$, the
surfaces limit to $K_1$ and as $x$ approaches $1$, the surfaces limit onto
$K_2$.  This family of surfaces is called a \textit{sweep-out}.  For each $x$,
the surface $\Sigma_x$ separates $M$ into a compression body $H^x_1$
and a handlebody $H^x_2$.

The fibers of the bundle structure on $M$ form a continuous one-parameter family of embedded, pairwise disjoint essential surfaces $\{T_y : y \in S^1\}$. Assume
the surfaces $\{\Sigma_x\}$ and $\{T_y\}$ are in general position.
The {\em Rubinstein-Scharlemann graphic} is the subset $R$ of
$(-1,1) \times S^1$ consisting of pairs $(x,y)$ such that $\Sigma_x$ and $T_y$
are tangent at some point in $M$.  General position implies that this set will be a
graph.

At each point $(x,y)$ in the complement of the graphic, the corresponding
surfaces $\Sigma_x$ and $T_y$ are transverse.  Label each point $(x,y)$ (in $(-1,1) \times S^1$) with a $1$ if some loop component of $\Sigma_x \cap T_y$ is essential in $\Sigma_x$ and is the boundary component of a disk or an essential annulus in $H^x_1$.  Note that for such an essential annulus, its other boundary component is a curve on $\bdry M$ isotopic to $\bdry T_y$.  Label $(x,y)$ with a $2$ if a loop is essential in $\Sigma_x$ and bounds a disk in $H^x_2$.  Given a second point $(x',y')$ in the same component of the complement of the graphic as $(x,y)$, a piecewise vertical an horizontal path from $(x,y)$ to $(x',y')$ determines an ambient isotopy of $M$ taking $\Sigma_x$ to $\Sigma_{x'}$ and $T_y$ to $T_{y'}$. Thus any two points in the same component of the complement of the graphic have the same labels.

\begin{claim} \label{claim:distinctlabelheight}
For a fixed $x$, there cannot be values $y,y'$ such that $(x,y)$ is labeled
with a $1$ and $(x,y')$ is labeled with a $2$. (In particular, a point cannot have
both labels $1$ and $2$.)
\end{claim}

\begin{proof}
Let $M'$ be the result of Dehn filling $\partial M$ along the slope of the boundary of a level surface $T_y$.  Then
$M'$ is a closed torus bundle and the image of a loop which is
boundary parallel in $T_y$ will be trivial in $M'$.  The
induced Heegaard splitting of $M'$ will be irreducible because a
torus bundle cannot be a lens space.  Moreover, an irreducible,
genus-two Heegaard splitting is strongly irreducible so the induced
Heegaard splitting of $M'$ is strongly irreducible.

Now assume for contradiction the point $(x,y)$ is labeled with a $1$ and $(x,y')$ is labeled with a $2$.  Then there is a loop in $\Sigma_x \cap T_y$ that bounds an essential disk in the filling of $H_1^x$ (after the filling, an essential annulus is capped off and becomes a disk) and a loop $\Sigma_x \cap T_{y'}$ that bounds an essential disk in $H_2^x$.  Because $T_y$ and $T_{y'}$ are disjoint, each loop in $\Sigma_x \cap T_{y'}$ is disjoint from each loop in $\Sigma_x \cap T_y$.  Thus the Heegaard surface $\Sigma_x$ for $M'$ is weakly reducible and therefore
reducible.  This contradicts the assumption that $M'$ is not a lens space.
\end{proof}

\begin{claim}\label{claim:TynotinH1}
No fiber $T_y$ can be made disjoint from $H^x_2$.
\end{claim}

\begin{proof}
If a fiber $T_y$ were disjoint from $H^x_2$, then it would be
contained in $H^x_1$.  Since $T_y$ is essential, a non-separating
compressing disk of $H^x_1$ must be disjoint from $T_y$.  Therefore
compressing $H^x_1$ along such a disk yields a manifold which is
homeomorphic to $T^2 \times I$ and contains a properly embedded,
essential, once-punctured torus, which cannot occur.
\end{proof}

\begin{claim}\label{claim:nonemptylabels}
In a component of $(-1,1) \times S^1 - R$ that intersects $(-1, -1+\epsilon) \times S^1$ for suitably small $\epsilon$, each point is labeled with a $1$.  In a component of $(-1,1) \times S^1 - R$ that intersects $(1-\epsilon, 1) \times S^1$ for suitably small $\epsilon$, each point is labeled with a $2$.
\end{claim}

\begin{proof}
Fix $y$.  Because $\Sigma_x$ limits onto $K_1$, if $x$ is near $-1$,
then $T_y \cap H^x_1$ will be a collection of disks and an essential
annulus.   Hence $(x,y)$ has label $1$.

If $x$ is near $1$, then $T_y \cap H^x_2$ consists of disks. If 
none of these disks are essential in $H^x_2$,
then $T_y$ can be pushed into $H_1^x$,
contradicting Claim~\ref{claim:TynotinH1}.  Hence $(x,y)$ has label
$2$.
\end{proof}

\begin{claim}\label{claim:unlabeledpoints}
For some value $x$, the point $(x,y)$ is unlabeled for every $y$.
\end{claim}

\begin{proof} 
The complement of the graphic is an open set, as is each component of the
complement.  The union of all the components labeled with a 1 projects to
an open set in $(-1,1)$ as does the union of the components labeled with a 2.
By Claim~\ref{claim:nonemptylabels} the projection of each set is non-empty.  By Claim~\ref{claim:distinctlabelheight} the images of these projections are disjoint.  Because $(-1,1)$ is connected, it cannot be written as the union of two non-empty open sets.  Thus there is a point $x$ in the complement of the projections.  This $x$ has the desired property.
\end{proof}

Fix an $x$ such that for every $y$, $(x,y)$ is unlabeled as
guaranteed by Claim~\ref{claim:unlabeledpoints}.  If for some $y$,
$T_y$ is transverse to $\Sigma_x$ and each loop of $\Sigma_x \cap
T_y$ is trivial in $T_y$ then $T_y$ can be isotoped into $H^x_1$. As
noted in Claim~\ref{claim:TynotinH1} above, this is impossible.
Therefore for each $y \in S^1$, if $T_y$ and $\Sigma_x$ intersect
transversely, then the intersection $\Sigma_x \cap T_y$ must contain
a loop which is essential in $T_y$.   Since two disjoint, essential
loops in a once-punctured torus are parallel, the essential loops in
$\Sigma_x \cap T_y$ are all parallel in $T_y$, for each $y$.

For each point $p \in \Sigma_x$ there is a $y \in S^1$ such that $p$ is
contained in $T_y$.  Define the function $f_x : \Sigma_x \rightarrow S^1$
such that $p$ is contained in $T_{f(p)}$ for each $p \in \Sigma_x$.  By general
position, either $f_x$ is a circle-valued Morse function or $f_x$ is a near-Morse function such that the critical points consist of either a single isolated degenerate critical point or just two critical points at the same level.

Assume for contradiction $f_x$ is a Morse function.  Then for each
$y \in S^1$, the level set $f_x^{-1}(y)$ contains a loop which is
essential in $\Sigma_x$.  Thus $\Sigma_x \cap T_y$ contains a loop
which is essential in $\Sigma_x$ and therefore essential in $T_y$.
By continuity, the isotopy class in $\Sigma_x$ of the essential
loops cannot change as $y$ varies, so the monodromy $\phi$ must send
this essential loop onto itself.  A quick calculation shows the
homology of this manifold would have rank 3, contradicting the fact
that $M$ admits a genus two Heegaard splitting.

If $f_x$ has an isolated degenerate critical point, then once again
$f_x^{-1}(y)$ must contain an essential loop for each $y$.  This leads to the
same contradiction as when $f_x$ is a Morse function, see~\cite{cooper-scharlemann}.  Therefore we conclude that $f_x$ must have two
critical points on the same level.  As described in~\cite{cooper-scharlemann},
this implies that $\Sigma_x$ is embedded in $M$ as shown in
Figure~\ref{fig:twocrits}.
\begin{figure}[htb]
\centering
\includegraphics[width=3.5in]{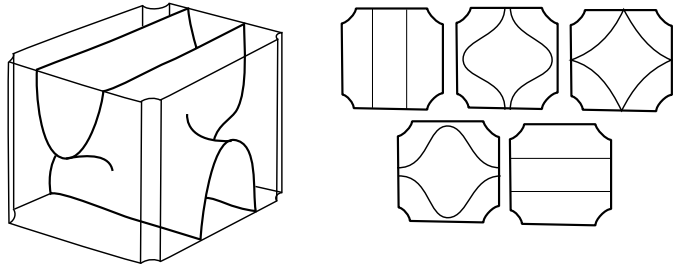}
\caption{$\Sigma_x$ in the once-punctured torus bundle $M$ cut open
along a fiber $T_y$, and the sequence of intersections of $\Sigma_x$
with the fibers $T_y$ as $y$ passes through the critical value of
$f_x$.} \label{fig:twocrits}
\end{figure}

The right side of the figure shows how the loops of $\Sigma_x \cap T_y$ sit in each $T_y$.  As $y$ passes through the critical value $c$ of $f_x$ from $c+\epsilon$ to $c-\epsilon$, the isotopy class $\ell_{c+\epsilon}$ of the essential loops of $\Sigma_x \cap T_{c+\epsilon}$ is replaced by a new isotopy class $\ell_{c-\epsilon}$ of the loops of $\Sigma_x \cap T_{c-\epsilon}$ which, under the projection $T \times (c-\epsilon, c+\epsilon) \to T$, intersects the original in a single point.  Because this is the only level at which this happens, we conclude that there are essential
loops $\alpha_1, \alpha_2 \subset T$ such that $\alpha_1 \cap \alpha_2$ is a
single point and $\phi(\alpha_1) = \alpha_2$.

Let $s_1$ be a Dehn twist along $\alpha_1$ and $s_2$ a Dehn twist
along $\alpha_2$.  Either composition $(s_2 s_1 s_2)^{\pm 1}$ takes $\alpha_1$ to $\alpha_2$ and $\alpha_2$ to $\alpha_1$.  By choosing $+1$ or $-1$ appropriately, we can ensure that the map sends $\alpha_1$ to $\alpha_2$ with the same orientation as $\phi$.  Composing further with $s_1^n$ for some $n$ will take $\alpha_2$ to $\phi(\alpha_2)$.  Thus $(s_2 s_1 s_2)^{\pm 1}s_1^n$ is isotopic to $\phi$.
\end{proof}

 \section{Genus one fibered knots via closed $3$--braids}\label{sec:3}

A genus one fibered knot may be described in terms of a double
branched covering of a closed 3-braid whose word encodes the
monodromy of the fibering. This viewpoint allows us to describe
tunnel number one once-punctured torus bundles.

\begin{lemma}\label{lem:corresp}
Every genus one fibered knot is the image of the braid axis of a
closed $3$--braid in $S^3$ in the double branched cover of the
closed $3$--braid.  Conversely, the image of the braid axis of a closed
$3$--braid in the double branched cover of the closed $3$--braid is
always a genus one fibered knot.
\end{lemma}

By considering the standard involution of the once-punctured torus
$T$, we will show that the set of genus one fibered knots in
$3$--manifolds and the set of braid axes of closed $3$--braids in $S^3$ are in
one-to-one correspondence.  We sketch the passage from genus one fibered knots
to axes of closed $3$--braids in $S^3$ below.  The return passage is
then clear.  (See Section 5 of \cite{birman-hilden:hsobcoS3},
Sections 4 and 5 of \cite{mr:gfm}, and Section 2 of
\cite{baker:cgofkils}.)

\begin{proof}[Sketch]
The standard involution $\iota$ with three fixed points of the
once-punctured torus $T$ extends to an involution of the
once-punctured torus bundle $M$ that takes a fiber to a fiber. Since
$\iota$ commutes with the monodromy of $M$, quotienting $M$ by this
involution yields a closed $3$--braid $\hat{\beta}$ in a solid torus
$W$ where the braid $\hat{\beta}$ is the image of the fixed set of
$\iota$ and a meridional disk of $W$ is the image of a fiber.

The extension of $\iota$ is a fixed point free involution on the
boundary of $M$, so $\iota$ further extends to an involution of the
filling $M'$ of $M$ by a solid torus $V'$ whose meridians intersect
fibers of $M$ just once.  In $M'$ the core of $V'$ is a genus one
fibered knot $K$.  Since the extension of $\iota$ acts as a free
involution of $V'$, in the quotient $V'$ descends to a solid torus
$V$.  The meridian of $V$ intersects the meridian of $W$ once, as
the meridian of $V'$ intersects the boundary of a fiber of $M$ once.
Hence under this quotient, $K$ descends to the core of $V$ which is
the axis of the closed $3$--braid $\hat{\beta}$ in $V \cup W \cong
S^3$.  This defines the correspondence.
\end{proof}

\begin{lemma}\label{lem:tnobraid}
Every tunnel number one, genus one fibered knot is the lift of the
braid axis of the closure $\hat{\beta}_{k,n}$ of the braid
$\beta_{k,n} = (\sigma_2 \sigma_1 \sigma_2)^k \sigma_1^n$ in the
double cover of $S^3$ branched over $\hat{\beta}_{k,n}$, where $k$
is odd and $n$ is an arbitrary integer.
\end{lemma}

\begin{figure}
\centering
\input{betakn.pstex_t}
\caption{The braid $\beta_{k,n} = (\sigma_2 \sigma_1 \sigma_2)^k \sigma_1^n$.}
\label{fig:betakn}
\end{figure}
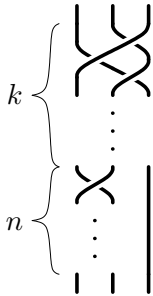

A depiction of the braid $\beta_{k,n} = (\sigma_2 \sigma_1 \sigma_2)^k \sigma_1^n$ is shown in Figure~\ref{fig:betakn}.

\begin{proof}
Let $\widehat{M}$ be a $3$--manifold that contains a tunnel number
one, GOF-knot $K$. Let $M$ be the once-punctured torus bundle
exterior $\widehat{M}-N(K)$ with fiber $T$.

By Lemma~\ref{monodromylem}, the monodromy of $M$ may be given by
$(s_2 s_1 s_2)^{\pm 1}s_1^n$ up to Dehn twists along the boundary of
a fiber $T$.  Write a single Dehn twist along $\bdry T$ as $(s_2 s_1
s_2)^{4}$.  Then the monodromy $T \to T$ fixing $\bdry T$ with an
orbit on $\bdry M$ bounding a meridional curve for $K$ is $(s_2 s_1
s_2)^{4\ell \pm 1}s_1^n$ for some $\ell \in \Z$.

By Lemma~\ref{lem:corresp}, $\widehat{M}$ is the double branched
cover of a closed $3$--braid in $S^3$ with braid axis lifting to
$K$.  Accordingly, $M$ is the double branched cover of the closed
$3$--braid in the solid torus exterior of the braid axis.  A
meridian of $K$ in $\bdry M$ corresponds to a meridian of the braid
axis in the solid torus.

Under the quotient of the covering involution, the Dehn twist $s_i$
along $\alpha_i$ corresponds to the braid $\sigma_i$, a right-handed
crossing between the $i$th and $(i+1)$th strands.  Thus $K$ in
$\widehat{M}$ corresponds to the braid axis of the closed braid
$(\sigma_2 \sigma_1 \sigma_2)^{k} \sigma_1^n$ where $k=4\ell \pm1$.
\end{proof}

\begin{remark}\label{rem:dbcsurgery}
A meridian of the braid axis is a longitudinal curve on the solid
torus containing the closed $3$--braid.  In the double branched
cover it lifts to two meridians of the genus one fibered knot.  The
longitude of the braid axis is the meridian of the solid torus
containing the closed $3$--braid.  It lifts to the longitude of the
genus one fibered knot.

In these coordinates for the braid axis, a slope of $p/q$ lifts to
the slope $2p/q$.  If $q$ is even (and $p$ and $q$ are coprime) this
slope is to be interpreted as two parallel curves of slope
$p/(q/2)$.  It follows that $1/q$ surgery on a genus one fibered
knot corresponds to inserting $2q$ full twists (right-handed if
$q<0$, left-handed if $q>0$) into the $3$--braid.

If $q$ is odd then a $p/q$ slope lifts to a single curve of slope
$2p/q$.  Hence $p/q$ surgery on the braid axis cannot lift to
surgery on the genus one fibered knot in the double branched cover
unless the core of the surgery solid torus is added to the branch
locus.
\end{remark}

\section{Genus one fibered knots in lens spaces}\label{sec:4}

By Corollary 4.12 of \cite{hodgson-rubinstein:iaiols}, the lens
space $L(\alpha, \beta)$ is the double cover of $S^3$ branched over
a link $L$ if and only if $L$ is equivalent to the two bridge link
$\B(\alpha, \beta)$.  Thus to understand GOF-knots in lens spaces
using Lemma~\ref{lem:corresp} we must consider representations of
two bridge links as closed $3$--braids; see
\cite{baker:cgofkils}.

Murasugi and later Stoimenow describe which two bridge
links admit such closed braid representations.  These two descriptions may be seen to be equivalent by working out their corresponding continued fractions.  Let ${\bf b}(L)$ denote the braid index of the link $L$.

\begin{prop}[Proposition 7.2, \cite{murasugi:otbioal}] \label{prop:murasugi}
Let $L$ be a two bridge link of type $\B(\alpha, \beta)$, where
$0<\beta<\alpha$ and $\beta$ is odd.  Then
\begin{enumerate}
\item ${\bf b}(L) = 2$ iff $\beta = 1$.
\item ${\bf b}(L) = 3$ iff either
\begin{enumerate}
\item \label{item:one} for some $p,q >1$, $\alpha = 2pq+p+q$ and $\beta = 2q+1$, or
\item \label{item:two} for some $q>0$, $\alpha =2pq+p+q+1$ and $\beta = 2q+1$.
\end{enumerate}
\end{enumerate}
\end{prop}

\begin{lemma}[Corollary 8, \cite{stoimenow:tspoc3b}]\label{lemma:stoimenow}
 If $L$ is a two bridge link of braid index 3, then $L$ has Conway notation $(p,1,1,q)$ or $(p,2,q)$ for some $p,q>0$.
 \end{lemma}

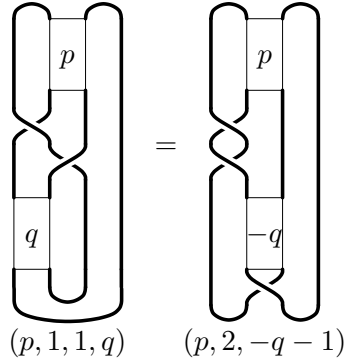
\begin{figure}
\centering
\input{twobridgeequiv.pstex_t}
\caption{The two bridge link $(p,1,1,q)$ is equivalent to $(p,2,-q-1)$.} \label{fig:twobridgeequiv}
\end{figure}

\begin{figure}
\centering
\input{p2qknot.pstex_t}
\caption{Closed $3$--braid representatives of the two bridge
link $(p,2,q)$.} \label{fig:twobridgethreebraid}
\end{figure}
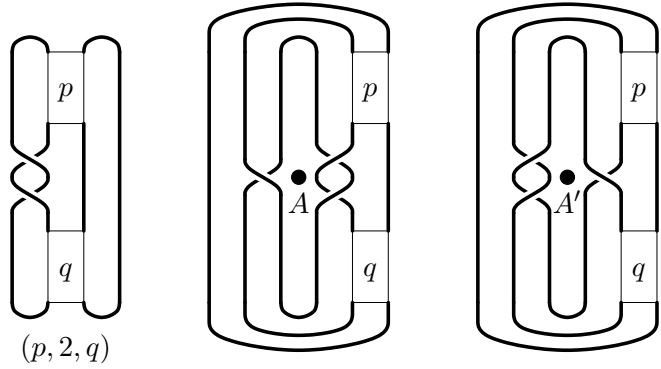

The two bridge link with Conway notation $(p,1,1,q)$ is shown in Figure~\ref{fig:twobridgeequiv}.  The link with Conway notation $(p,2,q)$ is shown in Figure~\ref{fig:twobridgethreebraid}.

In the other direction, we can determine which $3$--braids have closures that are two bridge links.

\begin{lemma}\label{lem:braidconj}
The closure of a $3$--braid $\beta$ is a two bridge link if
and only if $\beta$ or its mirror image is conjugate to $\sigma_2^{-1} \sigma_1^p \sigma_2^2 \sigma_1^q$ for some $p, q \in \Z$.
\end{lemma}

\begin{proof}
Assume the closure $\hat{\beta}$ of a $3$--braid $\beta$ is a two bridge link.  By Lemma~\ref{lemma:stoimenow} $\hat{\beta}$ may be described with Conway notation $(p,1,1,q)$ or $(p,2,q)$ for some integers $p,q$.  Figure~\ref{fig:twobridgeequiv} shows that a two bridge link with Conway notation $(p,1,1,q)$ is equivalent to one with Conway notation $(p,2,-q-1)$.  Figure~\ref{fig:twobridgethreebraid} shows the expression of a two bridge link with Conway notation $(p,2,q)$ as the closure of  the braid $\sigma_2^{-1} \sigma_1^p \sigma_2^2 \sigma_1^q$.  By Theorem~2.4 in~\cite{baker:cgofkils} the braid axis is unique in Case (2) of Proposition~\ref{prop:murasugi}.  For Case (1), the braid axes are classified and have the desired form (see the proof of Theorem~2.4 in~\cite{baker:cgofkils}).

If a $3$--braid $\beta$ is conjugate to $\sigma_2^{-1} \sigma_1^p
\sigma_2^2 \sigma_1^q$ for some $p, q \in \Z$, then its closure can readily be identified with a two bridge link, cf.\ Figure~\ref{fig:twobridgethreebraid}.
\end{proof}

\comment{
\begin{remark}
Figure~\ref{fig:twobridgethreebraid} shows the two bridge link with Conway notation $(p,2,q)$ can be expressed as closures of the
(often non-conjugate) braids $\sigma_2^{-1} \sigma_1^p \sigma_2^2
\sigma_1^q$ and $\sigma_2^{-1} \sigma_1^q \sigma_2^2 \sigma_1^p$.
\end{remark}
}

\begin{lemma}\label{lem:3braidlensspace}
The double branched cover of the closure of the braid $\sigma_2^{-1}
\sigma_1^p \sigma_2^2 \sigma_1^q$ is the lens space
$L(2pq+p+q,2q+1)$.
\end{lemma}
\begin{proof}
Since the closure of the braid $\sigma_2^{-1} \sigma_1^p \sigma_2^2
\sigma_1^q$ is a two bridge link with Conway notation $(p,2,q)$ (see Figure~\ref{fig:twobridgethreebraid}), it
corresponds to the two bridge link $\B(\alpha, \beta)$ where
$\alpha/\beta = p+\cfrac{1}{2+\frac{1}{q}}=\frac{2pq+p+q}{2q+1}$.
Since the double branched cover of the two bridge link $\B(\alpha,
\beta)$ is the lens space $L(\alpha, \beta)$ the result follows.
\end{proof}

\section{Tunnel number one, genus one fibered knots in lens spaces}\label{sec:5}

The monodromy of tunnel number one once-punctured torus bundle has a
special form, which can now be connected with lens spaces via
two bridge links.

\begin{lemma}\label{lem:braidform}
The closure $\hat{\beta}_{k,n}$ of the braid $\beta_{k,n} =
(\sigma_2 \sigma_1 \sigma_2)^k \sigma_1^n$, where $k$ is odd, is a
two bridge link if and only if $k=\pm1$, $(k,n) = \pm(-3,3)$, or $(k,n)=\pm(-3,5)$.
\end{lemma}

\begin{proof}
If $k = \pm1$ then
\[\beta_{\pm1,n} = (\sigma_2 \sigma_1 \sigma_2)^{\pm1} \sigma_1^n = (\sigma_1 \sigma_2 \sigma_1)^{\pm1} \sigma_1^n \equiv \sigma_2^{\pm1} \sigma_1^{n\pm2}.\]
Thus $\hat{\beta}_{\pm1,n}$ is the two bridge link $\B(n\pm2, 1)$
(the $(2,n\pm2)$--torus link).  The reader may check that $\hat{\beta}_{k,n}$ is the two bridge knot $\B(\mp5,1)$ for $(k,n) = \pm(-3,3)$ and the two bridge knot $\B(\pm7,2)$ for $(k,n)=\pm(-3,5)$.

If $n$ is even then $\hat{\beta}_{k,n}$ has one unknotted component
and one potentially knotted component (a $(2,k)$--torus knot).
Since a two bridge link has either one component or two unknotted
components, the $(2,k)$--torus knot must be the unknot.  Hence $k =
\pm 1$.

Assume that $n$ is an odd integer.  By taking mirror images, we
only need consider the case that $k \equiv 1 \mod 4$.  We will
deal with the three cases $n=1, n=-1,$ and $|n|>2$ separately.

Since $k = 4 \ell + 1$ for some $\ell \in \Z$ we have $\beta_{k,n} =
(\sigma_2 \sigma_1 \sigma_2)^{4\ell} \beta_{1,n}$.  By
Remark~\ref{rem:dbcsurgery}, the double branched cover of
$\hat{\beta}_{k,n}$ may be obtained from the double branched cover
$L(n+2,1)$ of $\hat{\beta}_{1,n}$ by $(1/\ell)$--Dehn surgery on the
lift of its braid axis.

If $n=-1$, then the braid axis of $\hat{\beta}_{1,n}$ lifts to the
right handed trefoil. A $(1/\ell)$--Dehn surgery on this knot this
is a lens space only when $\ell=0$ \cite{moser:esaatk}; hence $k=1$.

To handle the cases when $n=+1$ and $|n|>2$, we must further
understand the relationship between the order of the lens space and
the braid structure of $\hat{\beta}_{k,n}$.
The exponent sum of a braid is the sum of the exponents of a word in the letters $\{\sigma_1, \sigma_2\}$ representing the braid.
By Lemma~\ref{lem:braidconj}, $\hat{\beta}_{k,n}$ is a two bridge
link if and only if it is conjugate to $\sigma_2^{-1} \sigma_1^p
\sigma_2^2 \sigma_1^q$.  Since an exponent sum is invariant under conjugation, we must have:
\begin{equation}\label{eqn:1}
3k+n = p+q+1.
\end{equation}
Since the double branched cover of $\hat{\beta}_{4\ell+1,n}$ is obtained from
$L(n+2,1)$ by $(1/\ell)$--Dehn surgery on the lift of the braid axis, its first homology group must be cyclic of order $n+2$.
By Lemma~\ref{lem:3braidlensspace}, the double branched cover of
closure of the braid $\sigma_2^{-1} \sigma_1^p \sigma_2^2
\sigma_1^q$ is the lens space $L(2pq+p+q,2q+1)$.  Thus the orders of
first homology of the double branched cover of each
$\hat{\beta}_{4\ell+1,n}$ and the closure of $\sigma_2^{-1}
\sigma_1^p \sigma_2^2 \sigma_1^q$ must agree:
\begin{equation}\label{eqn:2}
|n+2| = |2pq+p+q|.
\end{equation}
Since $n$ is assumed to be odd, Equation~(\ref{eqn:2}) implies that
the integers $p$ and $q$ have different parity.

Putting Equations~(\ref{eqn:1}) and (\ref{eqn:2}) together, we
obtain $|p+q+3-3k| = |2pq+p+q|$.  Hence either
\begin{equation}\label{eqn:3}
3k-3 = 2(pq+p+q) \mbox{ or } 3k-3 = -2pq.
\end{equation}

Now consider the case when $n=+1$.
For $\hat{\beta}_{4\ell+1,1}$ to be a two bridge knot,
Equation~(\ref{eqn:2}) implies that $\{p,q\}=\{0,3\}, \{0,-3\},
\{-1, 2\}$, or  $\{-1, -4\}$.  Since $n=+1$, Equation~(\ref{eqn:1})
then implies that $k = 1, -1, 1/3$, or $-5/3$ respectively.  Because
$k$ must be an integer, $k = \pm1$.

Lastly consider the case when $|n| > 2$.  Denote the lift of the
braid axis of $\hat{\beta}_{1,n}$ by $K_{1,n}$. The action of the
monodromy on the homology of the fiber of $K_{1,n}$ is given by
\[ \left(\begin{array}{cc} 1 & 0 \\ -1 & 1\end{array} \right)
\left(\begin{array}{cc} 1 & n+2 \\ 0 & 1\end{array} \right) =
\left(\begin{array}{cc} 1 & n+2 \\ -1 & -n-1\end{array} \right). \]
This has trace $-n$.  By \cite{cb:aosanat} this monodromy is
pseudo-Anosov  for $|n|>2$ and therefore by \cite{Th2} the
complement of $K_{1,n}$ is hyperbolic. According to the Cyclic
Surgery Theorem \cite{cgls:dsok}, $1/\ell$ surgery on $K_{1,n}$ is a
lens space only if $\ell = 0$, $\ell = +1$, or $\ell = -1$.  Thus
only for $k=1$, $k=5$ or $k=-3$ respectively could
$\hat{\beta}_{k,n}$ with $|n|
> 2$ be a two bridge link.  We complete the proof with an examination of the latter two possibilities.

{\bf Case A: $k=5$}

By Equation~(\ref{eqn:3}) either $-6 = pq$ or $6 = pq+p+q$.  Thus
$-6=pq$ or $7=(p+1)(q+1)$.  Since $p$ and $q$ have opposite parity,
$-6=pq$, and so $\{p, q\} = \{-2, 3\}, \{2, -3\}, \{1, -6\},$ or $\{-1,
6\}$. Hence by Equation~(\ref{eqn:1}) $n=-13, -15, -19,$ or  $-9$ respectively.  {\tt Braid Group Calculator} \cite{bgc} shows that none of the braids $\beta_{5,n}$ are conjugate to $\sigma_2^{-1} \sigma_1^{p} \sigma_2^2 \sigma_1^{q}$ for corresponding choices of $n$ and $\{p,q\}$.  By Lemma~\ref{lem:braidconj} this implies
the braid closure $\hat{\beta}_{5, n}$ with $|n|>2$ is not a two
bridge link.

{\bf Case B: $k=-3$}

By Equation~(\ref{eqn:3}) either $6=pq$ or $-6=pq+p+q$.  Thus
$6=pq$ or $-5=(p+1)(q+1)$.  Since $p$ and $q$ have opposite parity,
$6=pq$, and so $\{p, q\} = \{2,3\}, \{-2, -3\}, \{1, 6\}$, or $\{-1,
-6\}$. Hence by Equation~(\ref{eqn:1}) $n= 15, 5, 17$, or $3$ respectively.  {\tt Braid Group Calculator} \cite{bgc} shows that the braids $\beta_{-3,n}$ are not conjugate to $\sigma_2^{-1} \sigma_1^{p} \sigma_2^2 \sigma_1^{q}$ for $(n,\{p,q\}) = (15,\{2,3\})$ and $(n,\{p,q\}) = (17, \{1,6\})$.  However these braids are conjugate if $(n,\{p,q\}) = (5,\{-2,-3\})$ or $(n,\{p,q\}) = (3, \{-1,-6\})$.
By Lemma~\ref{lem:braidconj} this implies the braid closure $\hat{\beta}_{-3, n}$ with $|n|>2$ is not a two bridge link unless $n = 3$ or $n = 5$.
\end{proof}

\begin{lemma}\label{lem:exception}
The braid $\beta_{-3,5}$ is not conjugate to a braid $\beta_{\epsilon,n}$ for any choice of $\epsilon = +1$ or $-1$ and integer $n$;  neither is $\beta_{3,-5}$.
\end{lemma}

\begin{proof}
The proof of Lemma~\ref{lem:braidform} shows that $\beta_{-3,5}$ is conjugate to $\sigma_2^{-1} \sigma_1^{-2} \sigma_2^2
\sigma_1^{-3}$.  By Lemma~\ref{lem:3braidlensspace}, the double branched cover of its closure is the lens space $L(7,2)$.  For $\epsilon = \pm1$ the double branched cover of the closure of $\beta_{\epsilon,n}$ is the lens space $L(n+2\epsilon,1)$.  If $\beta_{-3,5}$ were conjugate to $\beta_{\epsilon,n}$ then the double branched covers of their closures would be equal, but this is not the case.  Similarly, $\beta_{3,-5}$ is not conjugate to $\beta_{\epsilon,n}$ for any choice of $\epsilon$ and $n$.
\end{proof}

\begin{remark}\label{rem:conjbraids}
As one may check, the braids $\beta_{-3,3}$ and $\beta_{-1,-3}$ are conjugate as are $\beta_{3,-3}$ and $\beta_{1,3}$.
\end{remark}

\section{Proof of Theorem~\ref{thm:main}}\label{sec:6}

Using our understanding of the braid structure associated to the
knot, we can now determine the precise lens space in which the knot
lies, as well as describe the structure of the knot.

\begin{thmmain}
Every tunnel number one, genus one fibered knot in a lens space is
the plumbing of a $r$--Hopf band and a $(\pm 1)$--Hopf band in the
lens space $L(r,1)$ with one exception.  Up to mirroring, this exception is the genus one fibered knot in $L(7,2)$ that arise as $(-1)$--Dehn surgery on the plumbing of a $7$--Hopf band and a $(+1)$--Hopf band.
\end{thmmain}

\begin{proof}
Let $K$ be a tunnel number one, genus one fibered knot in a lens space.  By Lemma~\ref{lem:corresp} $K$ corresponds to the braid axis of the closure of a $3$--braid $\beta$.  As $K$ lies in a lens space, $\hat{\beta}$ must be a $2$--bridge link \cite{hodgson-rubinstein:iaiols}.  Since $K$ is tunnel number one, Lemma~\ref{lem:tnobraid}, Lemma~\ref{lem:braidform}, and Remark~\ref{rem:conjbraids} together imply that $\beta$ must be conjugate to $\beta_{k,n}$ for either $k = \pm1$ and $n \in \Z$ or $(k,n) = \pm(-3,5)$.

Let us first assume $\beta = \beta_{\pm1, n}$.  Then $\beta$ is conjugate to $\sigma_2^{\pm 1} \sigma_1^{n \pm 2}$.  Setting $r = n \pm 2$, the closure $\hat{\beta}$ is then the $(2, r)$--torus link and its double branched cover is the lens space $L(r,1)$.

We may now view $K$ as the plumbing of an $r$--Hopf band and a
$(\pm1)$--Hopf band as follows.  First observe that $\beta$ is a
Markov stabilization of the $2$--braid $\sigma_1^r$: a third
string and a single crossing ($\sigma_2^{\pm 1}$) is added to
$\sigma_1^r$ to form the $3$--braid $\beta$ preserving its
closure.

As in Proposition 12.3 of \cite{bz:knots}, the lens space $L(r,1)$ can be decomposed along a torus arising as the double branched cover of a bridge sphere $S^2$ for the closure $\widehat{\sigma_1^r}$.  Isotope the braid axis $h$ of $\sigma_1^r$ to lie on $S^2$, and let $D$ be a component of $S^2 - h$.

The link $\widehat{\sigma_1^r}$ punctures $D$ twice, and thus in the double branched covering, $D$ lifts to an annulus $A$ bounded by the preimage of $h$.  By Lemma~11.8 of \cite{bz:knots}, the braid $\sigma_1^r$ corresponds to $r$
Dehn twists about the core of $A$.  In other words, the preimage of $h$ lifts to an $r$--Hopf link.

As demonstrated in Theorem 5.3.1 of \cite{montesinosmorton:flfcb}, the Markov stabilization $\beta$ of $\sigma_1^r$ corresponds to the plumbing of a Hopf band onto $A$ in the double branched cover.  Thus we see $K$ as the boundary of the plumbing of an $r$--Hopf band and a $(\pm1)$--Hopf band in the lens space $L(r,1)$.

Now assume $\beta$ is conjugate to $\beta_{-3,5}$.  Since $\hat{\beta}_{-3,5} = \B(7,2)$, $K$ is a knot in $L(7,2)$.
By Lemma~\ref{lem:exception}, $\beta$ is not conjugate to $\beta_{\pm1,n}$.  Therefore $\beta$ is not a Markov stabilization of any $2$--braid, and so $K$ is not the boundary of the plumbing of an $r$--Hopf band and a $(\pm1)$--Hopf band.  (Indeed, $K$ is not the Murasugi sum of any two fibered links.)  Nevertheless, observe that adding two full positive twists to $\beta$ gives the braid $(\sigma_2 \sigma_1 \sigma_2)^4 \beta_{-3,5} = \beta_{1,5}$.  Remark~\ref{rem:dbcsurgery} thus implies that $K$ is the core of a $(-1)$--Dehn surgery on the lift of the braid axis of $\hat{\beta}_{1,5}$ to the double branched cover, i.e.\ the boundary of the plumbing of a $7$--Hopf band and a $(+1)$--Hopf band.
\end{proof}

\section{Punctured torus bundles as knot complements.}\label{sec:7}

\begin{main2}
Every tunnel number one once-punctured torus bundle is the complement of
a GOF-knot in a lens space of the form $L(r,1)$.
\end{main2}

\begin{proof}
Let $M$ be a once-punctured torus bundle with tunnel number one.  Filling
$\bdry M$ along a slope that intersects each fiber once produces a
closed $3$--manifold $\widehat{M}$ in which the core of the filling
solid torus is a tunnel number one GOF-knot $K$.  By
Lemma~\ref{lem:tnobraid}, $K$ is the lift of the braid axis of the
closure of the braid $(\sigma_2 \sigma_1 \sigma_2)^k \sigma_1^n$
where $k=4\ell \pm1$ is odd.  Following Remark~\ref{rem:dbcsurgery}, we
may perform $-1/\ell$ surgery on $K$ to produce a manifold
$\widehat{M}'$ in which the core of the surgery solid torus is a
tunnel number one GOF-knot $K'$.  Hence $K'$ is the lift of the
braid axis of the closure of the braid $(\sigma_2 \sigma_1
\sigma_2)^{\pm1} \sigma_1^n$.  This closed braid is equivalent to
the two bridge link $\B(r,1)$ where $r=n\pm2$.  Therefore
$\widehat{M}'$, the double cover of this link, is the lens space
$L(r,1)$,  and $M$ is the complement of $K'$.
\end{proof}

\begin{main3}
Every tunnel number one once-punctured torus bundle is the $(r/1)$--Dehn filling of a boundary component of the exterior of the Whitehead link for some integer $r$.
\end{main3}

\begin{proof}
By Theorem~\ref{thm:main2}, a tunnel number one once-punctured torus bundle is the complement of a GOF-knot $K$ in a lens space $L(r,1)$.  By Theorem~\ref{thm:main}, $K$ is the plumbing of a $r$--Hopf band and a $(\pm 1)$--Hopf band.

The $r$--Hopf band is obtained by $r/1$ surgery on a circle $C$ that links an annulus $A$ whose boundary is the unlink in $S^3$.  Plumbing a $(\pm1)$--Hopf band onto $A$ produces the unknot $U$ whose union with $C$ forms the Whitehead link.  The $r/1$ surgery on $C$ transforms $U$ into the plumbing of a $(\pm1)$--Hopf band onto the $r$--Hopf band in the lens space $L(r,1)$.
Thus the image in $L(r,1)$ of $U$ is our knot $K$.
Therefore the complement of $K$ in $L(r,1)$ is obtained from the complement of the Whitehead link $U \cup C$ by the filling of the boundary component coming from $C$.
\end{proof}

\bibliographystyle{abbrv}

\end{document}

%% file: betakn.pstex_t
\begin{picture}(0,0)%
\epsfig{file=betakn.pstex}%
\end{picture}%
\setlength{\unitlength}{2960sp}%
\begingroup\makeatletter\ifx\SetFigFont\undefined%
\gdef\SetFigFont#1#2#3#4#5{%
  \reset@font\fontsize{#1}{#2pt}%
  \fontfamily{#3}\fontseries{#4}\fontshape{#5}%
  \selectfont}%
\fi\endgroup%
\begin{picture}(1464,2466)(5512,-2644)
\put(5851,-1036){\makebox(0,0)[rb]{\smash{{\SetFigFont{11}{13.2}{\rmdefault}{\bfdefault}{\updefault}{\color[rgb]{0,0,0}$k$}%
}}}}
\put(5851,-2086){\makebox(0,0)[rb]{\smash{{\SetFigFont{11}{13.2}{\rmdefault}{\bfdefault}{\updefault}{\color[rgb]{0,0,0}$n$}%
}}}}
\end{picture}%

%% file: twobridgeequiv.pstex_t
\begin{picture}(0,0)%
\epsfig{file=twobridgeequiv.pstex}%
\end{picture}%
\setlength{\unitlength}{2960sp}%
\begingroup\makeatletter\ifx\SetFigFont\undefined%
\gdef\SetFigFont#1#2#3#4#5{%
  \reset@font\fontsize{#1}{#2pt}%
  \fontfamily{#3}\fontseries{#4}\fontshape{#5}%
  \selectfont}%
\fi\endgroup%
\begin{picture}(2698,2996)(4726,-4099)
\put(5251,-1636){\makebox(0,0)[b]{\smash{{\SetFigFont{11}{13.2}{\rmdefault}{\mddefault}{\updefault}{\color[rgb]{0,0,0}$p$}%
}}}}
\put(4951,-3136){\makebox(0,0)[b]{\smash{{\SetFigFont{11}{13.2}{\rmdefault}{\mddefault}{\updefault}{\color[rgb]{0,0,0}$q$}%
}}}}
\put(5251,-4036){\makebox(0,0)[b]{\smash{{\SetFigFont{11}{13.2}{\rmdefault}{\mddefault}{\updefault}{\color[rgb]{0,0,0}$(p,1,1,q)$}%
}}}}
\put(6901,-1636){\makebox(0,0)[b]{\smash{{\SetFigFont{11}{13.2}{\rmdefault}{\mddefault}{\updefault}{\color[rgb]{0,0,0}$p$}%
}}}}
\put(6901,-3136){\makebox(0,0)[b]{\smash{{\SetFigFont{11}{13.2}{\rmdefault}{\mddefault}{\updefault}{\color[rgb]{0,0,0}$-q$}%
}}}}
\put(6901,-4036){\makebox(0,0)[b]{\smash{{\SetFigFont{11}{13.2}{\rmdefault}{\mddefault}{\updefault}{\color[rgb]{0,0,0}$(p,2,-q-1)$}%
}}}}
\put(6076,-2386){\makebox(0,0)[b]{\smash{{\SetFigFont{11}{13.2}{\rmdefault}{\mddefault}{\updefault}{\color[rgb]{0,0,0}$=$}%
}}}}
\end{picture}%

%% file: p2qknot.pstex_t
\begin{picture}(0,0)%
\epsfig{file=p2qknot.pstex}%
\end{picture}%
\setlength{\unitlength}{2960sp}%
\begingroup\makeatletter\ifx\SetFigFont\undefined%
\gdef\SetFigFont#1#2#3#4#5{%
  \reset@font\fontsize{#1}{#2pt}%
  \fontfamily{#3}\fontseries{#4}\fontshape{#5}%
  \selectfont}%
\fi\endgroup%
\begin{picture}(5508,3046)(5851,-4024)
\put(8326,-2761){\makebox(0,0)[b]{\smash{{\SetFigFont{11}{13.2}{\rmdefault}{\mddefault}{\updefault}{\color[rgb]{0,0,0}$A$}%
}}}}
\put(8926,-1786){\makebox(0,0)[b]{\smash{{\SetFigFont{11}{13.2}{\rmdefault}{\mddefault}{\updefault}{\color[rgb]{0,0,0}$p$}%
}}}}
\put(8926,-3286){\makebox(0,0)[b]{\smash{{\SetFigFont{11}{13.2}{\rmdefault}{\mddefault}{\updefault}{\color[rgb]{0,0,0}$q$}%
}}}}
\put(6376,-3961){\makebox(0,0)[b]{\smash{{\SetFigFont{11}{13.2}{\rmdefault}{\mddefault}{\updefault}{\color[rgb]{0,0,0}$(p,2,q)$}%
}}}}
\put(6376,-1786){\makebox(0,0)[b]{\smash{{\SetFigFont{11}{13.2}{\rmdefault}{\mddefault}{\updefault}{\color[rgb]{0,0,0}$p$}%
}}}}
\put(6376,-3286){\makebox(0,0)[b]{\smash{{\SetFigFont{11}{13.2}{\rmdefault}{\mddefault}{\updefault}{\color[rgb]{0,0,0}$q$}%
}}}}
\put(10576,-2761){\makebox(0,0)[b]{\smash{{\SetFigFont{11}{13.2}{\rmdefault}{\mddefault}{\updefault}{\color[rgb]{0,0,0}$A'$}%
}}}}
\put(11176,-1786){\makebox(0,0)[b]{\smash{{\SetFigFont{11}{13.2}{\rmdefault}{\mddefault}{\updefault}{\color[rgb]{0,0,0}$p$}%
}}}}
\put(11176,-3286){\makebox(0,0)[b]{\smash{{\SetFigFont{11}{13.2}{\rmdefault}{\mddefault}{\updefault}{\color[rgb]{0,0,0}$q$}%
}}}}
\end{picture}%